\documentclass[conference]{IEEEtran}
\IEEEoverridecommandlockouts
\usepackage{cite}
\usepackage{amsmath,amssymb,amsfonts}
\usepackage[ruled,vlined]{algorithm2e}
\usepackage[noend]{algorithmic}
\usepackage{graphicx}
\usepackage{textcomp}
\usepackage{xcolor}
\usepackage{subcaption}
\usepackage{flushend}
\usepackage{balance}
\usepackage[a4paper, total={184mm,239mm}]{geometry}
\usepackage{tikz}
\usepackage[outdir=./Figure/]{epstopdf}
\def\BibTeX{{\rm B\kern-.05em{\sc i\kern-.025em b}\kern-.08em
    T\kern-.1667em\lower.7ex\hbox{E}\kern-.125emX}}

\DeclareMathOperator*{\argmin}{arg\,min}

\newcommand\copyrighttext{%
  \footnotesize $\copyright$ 2021 IEEE. Personal use of this material is permitted. Permission from IEEE must be obtained for all other uses, in any current or future media, including reprinting/republishing this material for advertising or promotional purposes, creating new collective works, for resale or redistribution to servers or lists, or reuse of any copyrighted component of this work in other works.}
\newcommand\copyrightnotice{%
\begin{tikzpicture}[remember picture,overlay]
\node[anchor=south,yshift=10pt] at (current page.south) {\fbox{\parbox{\dimexpr\textwidth-\fboxsep-\fboxrule\relax}{\copyrighttext}}};
\end{tikzpicture}%
}

\begin{document}
\title{Energy-Efficient Adaptive System Reconfiguration for Dynamic Deadlines in Autonomous Driving
}

\author{\IEEEauthorblockN{Saehanseul Yi\IEEEauthorrefmark{1},
Tae-Wook Kim\IEEEauthorrefmark{2}, Jong-Chan Kim\IEEEauthorrefmark{2}\IEEEauthorrefmark{3} and
Nikil Dutt\IEEEauthorrefmark{1}}
\IEEEauthorblockA{
\IEEEauthorrefmark{1}
Department of Computer Science,
University of California, Irvine, USA\\
}

\IEEEauthorblockA{
\IEEEauthorrefmark{2}
Graduate School of Automotive Engineering,
Kookmin University, Korea\\
}

\IEEEauthorblockA{
\IEEEauthorrefmark{3}
Department of Automobile and IT Convergence,
Kookmin University, Korea\\
\{saehansy, dutt\}@uci.edu,
\{dsd8135, jongchank\}@kookmin.ac.kr}
}

\maketitle

\begin{abstract}
The increasing computing demands of autonomous driving applications make energy optimizations critical for reducing battery capacity and vehicle weight.
Current energy optimization methods typically target traditional real-time systems with static deadlines, resulting in conservative energy savings that are unable to exploit additional energy optimizations due to dynamic deadlines arising from the vehicle's change in velocity and driving context.
We present an adaptive system optimization and reconfiguration approach
that 
dynamically 
adapts the scheduling parameters and processor speeds to satisfy dynamic deadlines while consuming as little energy as possible. 
Our experimental results with an autonomous driving task set from Bosch and real-world driving data
show energy reductions up to 46.4\% on average in typical dynamic driving scenarios compared with traditional static energy optimization methods, demonstrating 
great potential for dynamic energy optimization gains by exploiting dynamic deadlines.
\end{abstract}


\copyrightnotice
\section{Introduction}
\label{sec:intro}
Due to the enormous amount of computing required for autonomous vehicles, its inherent high energy consumption has become one of the major hurdles when designing such computing systems. For example, even a single vehicle computer with multiple central processing units (CPUs) and graphics processing units (GPUs) alone reportedly consumes more than 2~kW of power, reducing an electric vehicle's driving range up to 12\%~\cite{lin2018architectural}. With this challenge, some energy optimization methods have been proposed that at the same time try to satisfy the stringent real-time requirements of automotive systems~\cite{kehr2017parcus, xie2018reliability}. However, since they commonly assume just (fixed) {\em static deadlines}, they do not reflect recent autonomous driving applications with time-varying {\em dynamic deadlines}. Such applications include the localization system with its dynamic latency constraint as a function of velocity~\cite{reid2019localization} and the truck platooning system where its control response times are adjustable to varying driving conditions~\cite{davis2012stability}.

With the above motivation, we aim to develop an energy-efficient software optimization method and a runtime framework that can specifically exploit the dynamic nature of deadlines found in many  autonomous driving applications. This study specifically focuses on minimizing CPU energy consumption by using {\em dynamic voltage and frequency scaling} (DVFS). Although conventional automotive microcontrollers often lack such features, recent application processors for autonomous driving mostly support DVFS.

\begin{figure}
\centering
     \begin{subfigure}{0.5\textwidth}
         \centering
         \includegraphics[width=1\textwidth]{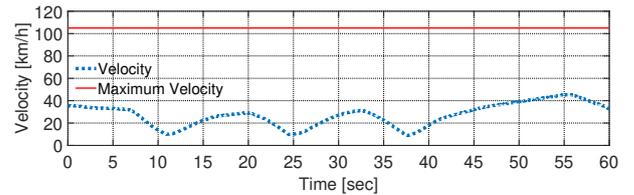}
         \caption{Example time history of vehicle velocity~\cite{velocity_dataset}}
         \vspace{0.3\baselineskip}
         \label{subfig-1:velocity)}
     \end{subfigure}
     \begin{subfigure}{0.5\textwidth}
         \centering
         \includegraphics[width=1\textwidth]{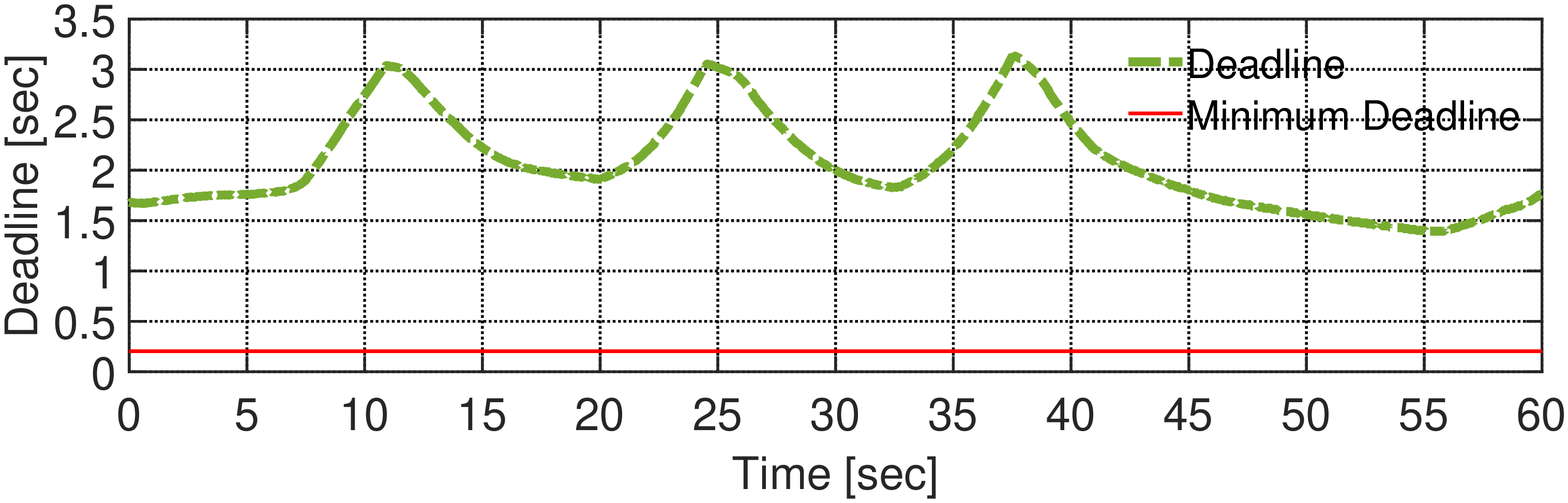}
         \caption{Corresponding dynamic deadlines ($\lambda$=20~m and $a^{max}$=2.5~ms\textsuperscript{-2})}
         \label{subfig-2:deadline}
     \end{subfigure}
        \caption{Example dynamic deadlines}
    \label{fig:velocity_deadline}
    \vspace{-0.55cm}
\end{figure}

To demonstrate our basic idea, Fig.~\ref{fig:velocity_deadline} shows an example time history of a vehicle's velocity and its corresponding dynamic deadline, assuming a velocity-deadline mapping function
\begin{equation}
d(v)=\frac{-v+\sqrt{v^2+2\lambda a^{max} }}{a^{max}},
\label{eq:dv}
\end{equation}
which represents the minimum time for a vehicle at its initial velocity $v$ to advance a fixed distance $\lambda$ assuming its maximum acceleration $a^{max}$. It is especially useful in truck platooning where trucks maintain a fixed longitudinal gap between them across various driving velocities~\cite{davila2011sartre} such that safe control decisions can be made more efficiently in terms of a maximum travel distance (i.e., $\lambda$) between sensing and actuation rather than by a rigid timing constraint. Hence, in the figure, the faster the vehicle runs, the shorter the deadline gets. The minimum deadline depicted by a red line is determined by the maximum velocity enforced by traffic regulations. Here, our basic idea is to trade the area between the time-varying dynamic deadline and the minimum deadline to reduce energy consumption by adaptively slowing down the CPU to the extent that guarantees the dynamic deadline, rather than adhering to the static deadline~\cite{kehr2017parcus, xie2018reliability} as in the traditional energy optimization methods.

In the automotive industry, complex control applications composed of multiple independent real-time tasks are commonly modeled with directed acyclic graphs (DAGs) whose nodes are periodic tasks and edges are read-write dependencies between tasks. Fig.~\ref{fig:Dag} shows an example DAG from Bosch~\cite{waters} for their reference autonomous driving system from sensors to actuators. Upon such a DAG, its dynamic deadlines are imposed by its worst-case end-to-end delays from sources to sinks, which is collectively decided by the periods of individual tasks. Then, our objective is to minimize the average power consumption while guaranteeing such dynamic deadlines. When doing that, we must satisfy (i) the system schedulability constraint (i.e., task periods) as well as (ii) the end-to-end deadline constraint. At first, we solve the problem by assuming a static deadline constraint. For that, we formulate it as a {\em geometric programming} (GP) problem~\cite{boyd2007tutorial}, which is a special form of non-convex optimizations that can be efficiently solved by a transformation into a convex problem.

To extend the above optimization method to time-varying dynamic deadlines, we partition the feasible deadline range into a number of discrete {\em modes}, where the system is separately optimized for each mode, assuming each mode's {\em shortest} deadline, respectively. Then, our runtime framework provides a safe mode change protocol that changes each task's period when the vehicle slows down or speeds up crossing across different modes. Our mode change protocol is designed not to miss any deadline if the mode change is from a shorter deadline to a longer one. However, we found that extra delays are unavoidable in the opposite direction (i.e., longer to shorter deadlines). Even in that case, however, we provide a mode change delay analysis method from which we can reserve enough safety margins to hide away the extra delays.

Our experimental results show that our approach reduces the average energy consumption up to 46.4\% in various real-world driving scenarios compared with the conventional method based on static deadlines. Moreover, our extensive simulation experienced no deadline miss due to our safe mode change protocol and delay analysis method. To the best of our knowledge, our work is one of the first attempts for optimizing energy consumption in computing systems for autonomous driving, explicitly focusing on dynamic deadlines.

This study's contributions can be summarized as follows:
\begin{itemize}
    \item We formulate an optimization problem for energy-efficient autonomous driving systems with time-varying dynamic deadlines and provide a GP-based optimal solution.
    \item We provide a safe mode change protocol that guarantees analyzable (if any) overheads, which can be safely manipulated in the design time.
\end{itemize}

The remainder of this paper is organized as follows: The next section describes the background and our problem. Section~\ref{sec:opt} presents our offline system optimization method. Section~\ref{sec:mode} explains the dynamic system reconfiguration approach. Section~\ref{sec:practical} discusses practical issues. Section~\ref{sec:experiments} provides the evaluation results. Section~\ref{sec:related} presents related work. Finally, Section~\ref{sec:conclusion} concludes the paper.

\section{Background and Problem Description}
\label{sec:back}

\subsection{System Model}
\label{sec:system}

This study assumes a computing system with a single CPU~\footnote{As an initial effort for the energy optimization based on dynamic deadlines, this study employs a rather simple system model, which will be extended in our future work. Its practical consideration will be discussed in Section~\ref{sec:practical}} supporting DVFS with a speed factor $s$ in the range of
\begin{equation}
0 < s_{min} \leq s\leq 1,
\end{equation}
where $s_{min}$ denotes the minimum speed factor used for the CPU idle time while 1 (i.e., 100\%) indicates the maximum processing speed. This study assumes that $s$ can be any continuous real value in that range such that general optimization techniques can be applied to. The practical consideration with commercial off-the-shelf (COTS) CPUs supporting only a number of discrete frequency levels will be discussed in Section~\ref{sec:practical}.

\begin{figure}
\centering
\includegraphics[width=9 cm]{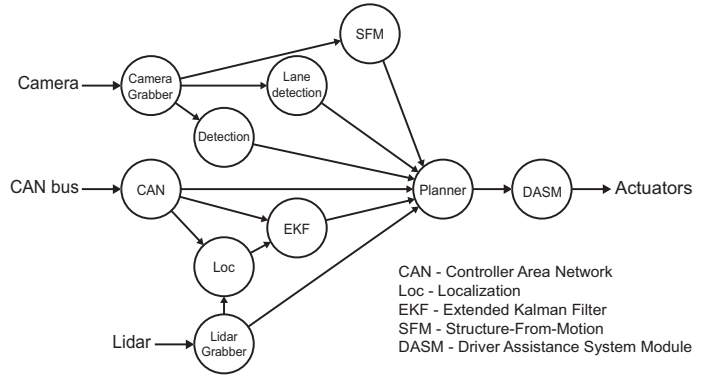}
\caption{A DAG from Bosch as a reference autonomous driving system  in the WATERS industrial challenge 2019~\cite{waters}}
\label{fig:Dag}
\vspace{-0.4cm}
\end{figure}

The system executes a set of $n$ implicit-deadline periodic tasks
\begin{equation}
V = \{\tau_1, \tau_2, \cdots, \tau_n\},
\end{equation}
where their read-write dependencies are represented by a DAG
\begin{equation}
G=(V, E \in V \times V)
\end{equation}
with its nodes $V$ and directed edges $E$ representing the task set and dependencies, respectively. Then, each task $\tau_i$ is characterized by its period $p_i$ and worst-case execution time (WCET) $e_i$ assuming $s$=1. Fig.~\ref{fig:Dag} shows an example DAG with ten tasks with complex dependencies beginning from sensors to actuators. Tasks communicate with each other with asynchronous message passing. Due to the multi-rate nature, oversampling or undersampling can happen to communication buffers, where newly arrived data always overwrite existing ones. This task model has been commonly used in many studies for automotive systems~\cite{feiertag2009compositional, davare2007period}.

When dealing with speed factors, we use the inter-task DVFS method, where $s$ can only be changed at context switching between tasks~\cite{bambagini2016energy}. Then, each task $\tau_i$ in the runtime is characterized by
\begin{equation}
\tau_i = (p_i, e_i, s_i),
\end{equation}
where $s_i$ is the per-task speed factor at the moment. Among them, only $e_i$ is a given value, whereas $p_i$ and $s_i$ are design variables. Thus, a complete {\em system configuration} $\pi$ can be described by a vector of tuples
\begin{equation}
\pi=((p_1, s_1), \cdots, (p_n, s_n)).
\end{equation}

We assume the earliest deadline first (EDF) scheduling algorithm such that the L\&L utilization bound~\cite{liu1973scheduling} can be used to test the system schedulability as
\begin{equation}
    U(\pi) = \sum_{i=1}^{n}\frac{e_i}{p_i s_i} \leq 100\%, \label{eq:ll_util}
\end{equation}
which has $s_i$ in the denominator reflecting the effective WCET that is inverse-linearly proportional to $s_i$. Based on the above schedulability test, we can guarantee every task's {\em periodicity}.

\subsection{Dynamic Deadlines}
\label{sec:dynamic}

In our system model, when referring to deadlines, they always mean the {\em end-to-end deadlines} from sensors to actuators, not the per-task implicit deadlines that are equal to periods. Thus, even when the system is schedulable, satisfying every task's $p_i$, it does not mean deadlines will be guaranteed. To formally define our notion of deadlines, let us assume $n_s$ {\em sensor tasks} (i.e., source nodes) and $n_a$ {\em actuator tasks} (i.e., sink nodes) in $G$. Then we say there are $n_s \times n_a$ unique {\em flows}, each of which has at least one {\em path} that is a sequence of adjacent tasks fully connecting a flow. The set of paths in $G$ is denoted by
\begin{equation}
\mathbb{P} = \{\delta_1, \delta_2, \cdots, \delta_{|\mathbb{P}|}\},
\end{equation}
where each path $\delta_i$ denotes an ordered set of task {\em indices} following the path. For example, Fig.~\ref{fig:Dag} has three $(3 \times 1)$ flows and eight paths. Then deadlines are imposed upon the paths such that newly arrived sensor data at time $t_1$ propagates through the DAG until it first gets out of an actuator task at time $t_2$ within a given deadline $d$ (i.e., $t_2 - t_1 \leq d$). In the automotive industry, the above notion is commonly referred to as {\em reaction time constraints}~\cite{extensions}.

\begin{figure}
\centering
\includegraphics[width=9 cm]{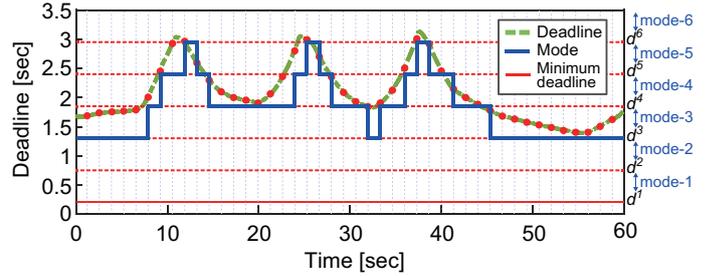}
\caption{Dynamic deadlines with discrete mode changes}
\label{fig:modes}
\vspace{-0.4cm}
\end{figure}

Fig.~\ref{fig:modes} shows continuous dynamic deadlines as the vehicle velocity changes, where vertical dashed lines depict discrete sensor arrivals. At each $k$-th sensor data arrival at time $t[k]$, its dynamic deadline $d[k]$ is decided as a red point by the vehicle velocity $v[k]$ with a given mapping function (e.g., $d(v)$ in Eq.~\eqref{eq:dv}). Thus, each sensor data arrival can be denoted by $(t[k], d[k]) \text{ for } k \geq 1$. Although many variables representing other physical states can be considered, this study focuses on the velocity-dependent deadlines as an initial effort toward a more general framework.

To efficiently manage dynamic deadlines, we employ a multi-mode approach, where a feasible deadline range is partitioned into $m$ discrete {\em modes}, where each mode guarantees the shortest deadline within its deadline range. For notational convenience, the modes are denoted by the per-mode shortest deadlines $\{d^1, d^2, \cdots, d^m\}$.

In Fig.~\ref{fig:modes}, its deadline domain is partitioned into six equal length modes, and at each sensor data arrival, the system mode is decided, possibly triggering mode changes. While the system is in a particular mode, the mode's shortest deadline is guaranteed, as depicted by the thick blue line.

\subsection{Power Model}
\label{sec:energy}
We use the following popular power model from~\cite{bambagini2016energy, bhuiyan2020energy}:
\begin{equation}
    P(s) = P_s + P_d(s) = \beta + \alpha s ^\gamma,
\label{eq:power_simple}
\end{equation}
where $P_s$ is the static power and $P_d(s)$ is the dynamic power parameterized by a speed factor $s$. The static power is expressed as $\beta$ independent of other parameters, whereas the dynamic power depends on $s$ while $\alpha$ and $\gamma \in [2, 3]$ are CPU-dependent parameters. In this work, we do not use CPU sleep states to reduce the static power, so our focus is to minimize the average dynamic power by minimizing $s$ as long as satisfying the dynamic deadlines.

For that, the average dynamic power can be calculated as follows: For each task $\tau_i$, its instantaneous dynamic power $\alpha s_i^\gamma$ is maintained during task's effective execution time $e_i/s_i$ affected by the speed factor. Since the power pattern repeats by its period $p_i$, the average power consumed in a unit time while executing $\tau_i$ can be calculated as
\begin{equation}
    P_i(p_i, s_i) = \frac{\alpha s_i^\gamma \times \frac{e_i}{s_i}}{p_i} = \frac{\alpha {s_i}^{\gamma-1} e_i}{p_i},
\end{equation}
which is a function of $p_i$ and $s_i$. Then, by summing up $n$ such average powers with the static power and the idle-time CPU power at $s_{min}$, the total average power for the system is given as a function of a system configuration $\pi$ as
\begin{equation}
\begin{aligned}
   P(\pi) = \beta + \alpha \sum_{i=1}^{n}\frac{{s_i^{\gamma-1}}{e_i}}{p_i} + \alpha s_{min}^\gamma\left(1 - \sum_{i=1}^{n}\frac{e_i}{p_i s_i}\right).
\end{aligned}
\label{eq:power}
\end{equation}

\subsection{Problem Description}
\label{sec:problem}

With the above system model, multi-mode deadline constraints, and energy model, our problem can be described as follows: Given a DAG $G$ of $n$ periodic tasks and $m$ discrete deadline constraints, our problem is to find the following two matrices:
\begin{equation}
\begin{pmatrix}
p_1^1 & p_2^1 & \cdots & p_n^1\\
\vdots & \vdots & \ddots&   \vdots\\
p_1^m & p_2^m &\cdots & p_n^m
\end{pmatrix}
\text{and}
\begin{pmatrix}
s_1^1 & s_2^1 & \cdots & s_n^1\\
\vdots & \vdots & \ddots&   \vdots\\
s_1^m & s_2^m &\cdots & s_n^m
\end{pmatrix},
\label{eq:matrices}
\end{equation}
where each $p_i^j$ and $s_i^j$ represent $\tau_i$'s optimal period and speed factor at the $j$-th mode, respectively, in terms of system-wide energy efficiency. Besides, the solution should satisfy {\em safe mode changes} such that the system can freely go back and forth between modes without violating the dynamic deadline requirements. For that, a safe runtime mode change protocol should be developed along with a mode change delay analysis method.

\section{Periods and Speed Factors Optimization}
\label{sec:opt}
This section formulates and solves the multi-mode system optimization problem. We begin by finding the optimal configuration assuming a single fixed mode without considering mode changes (Section~\ref{sec:per-mode}). Then, we extend the optimization method such that it can guarantee safe mode changes (Section~\ref{sec:across-modes}).

\subsection{Single-Mode Formulation}
\label{sec:per-mode}

This section explains how we can formulate a single-mode optimization as a baseline for the multi-mode system optimization.
As the objective function, the average power in Eq.~\eqref{eq:power} is used, without the constants $\alpha$ and $\beta$, where $\pi$ denotes two sets of decision variables  $p_i$s and $s_i$s with constrained domains as $p_i > 0$ and $s_{min} \leq s_i \leq 1$, respectively.

Then, we have two explicit constraints: (i) schedulability constraint as already discussed in Eq.~\eqref{eq:ll_util} and (ii) deadline constraint, for which we need to devise an end-to-end delay analysis model. We borrow a widely used model from~\cite{davare2007period}, which expresses the worst-case delay of a path $\delta$ as an accumulation of periods ($p_i$s) and worst-case response times (WCRTs) ($r_i$s) of every task in $\delta$, like the rightmost part in
\begin{equation}
\begin{aligned}
    D_\delta(\pi)=\sum_{i \in \delta}2p_i  \approx  \sum_{i \in \delta}(p_i + r_i),
\end{aligned}
\end{equation}
where $D_\delta(\pi)$ denotes the approximated worst-case delay of a path $\delta$ assuming a system configuration $\pi$. We approximate the original delay model to a linear form $\sum_{i \in \delta}2p_i$ as calculating the exact $r_i$ involves complex non-convex operations. Then, the deadline constraint $d$ is considered for every path $\delta$ in $\mathbb{P}$ (i.e., the set of paths in $G$).
Among the per-task delay components $2p_i$, one $p_i$ is for the waiting time until the task reads the sensor data ({\em waiting delay}), and the other is for processing the data ({\em processing delay)}.

Then our single-mode formulation is given as follows:
\begin{equation}
\begin{aligned}
    \underset{\pi}{\text{minimize}}
        & \quad P(\pi) = \sum_{i=1}^{n} {{s_i^{\gamma-1}} e_i \over {p_i}} + s^\gamma_{min}(1-\sum_{i=1}^{n} { {e_i}\over {p_i s_i}})  \\
    \text{subject to} 
        & \quad U(\pi) = \sum_{i=1}^{n}\frac{e_i}{p_i s_i} \leq 1 \\
        & \quad D_\delta(\pi) = \sum_{i\in\delta} 2p_i \leq d\, \ \ (\forall \delta \in \mathbb{P}).
\end{aligned}
\label{eq:single_mode}
\end{equation}

\subsection{Multi-Mode Formulation Considering Mode Changes}
\label{sec:across-modes}
Naively, the method in Section~\ref{sec:per-mode} can be repeatedly used to find every row of the multi-mode solution matrices in Eq.~\eqref{eq:matrices}. However, we cannot directly use this approach for a multi-mode system since it does not guarantee a safe transition between modes. 
Specifically, when old- and new-mode tasks coexist during a mode change, schedulability violations can occur even if each mode is schedulable in isolation~\cite{chen2018safemc}. 
Thus, we add a new constraint called \textit{per-task utilization invariability}, meaning every task's utilization $u_i=e_i/(p_is_i)$ is invariant across modes as
\begin{equation}
 \frac{e_i}{p_i^1 s_i^1}=\frac{e_i}{p_i^2 s_i^2}=\cdots=\frac{e_i}{p_i^m s_i^m}=u_i^* \text{ } (\forall i = 1,2,...,n),
\end{equation}
where $u_i^*$ denotes identical utilization for $\tau_i$ across modes. In that manner, even in the transient interval, the system's instantaneous utilization is maintained unchanged, which in turn guarantees the system schedulability~\cite{abdelzaher2004utilization}. Now, we use $u_i^*$s as our decision variables, replacing $s_i^j$s, which can be decided later by $s_i^j=e_i/(p_i^j u_i^*)$. Then our multi-mode optimization can be formulated as follows:
\begin{equation}
\begin{aligned}
    \underset{\hat{\pi}}{\text{minimize}}
        & \quad P(\hat{\pi}) = \sum_{j=1}^{m} { \sum_{i=1}^{n} {{e_i^\gamma} \over {(p_i^j)^\gamma (u_i^*)^{\gamma-1}}} + s^\gamma_{min}(1-\sum_{i=1}^{n} {u_i^*}}) \\
    \text{subject to} 
        & \quad U(\hat{\pi}) = \sum_{i=1}^{n} {u_i^*} \leq 1,  \\
        & \quad D^j_\delta(\hat{\pi}) = \sum_{i\in\delta} 2p_{i}^j \leq d^j\,\ \ (\forall j \in [1, m], \forall \delta \in \mathbb{P}),
\end{aligned}
\label{eq:across_mode}
\end{equation}
where $\hat{\pi}$ denotes the newly defined multi-mode system configuration with $p_i^j$s and $u_i^*$s. The objective function is the sum of average power in each mode, after eliminating the $\alpha$ and $\beta$ from Eq.~\eqref{eq:power} for the notational simplicity. The first constraint is the system schedulability now expressed by $u_i^*$s. The second constraint is for the dynamic deadlines across $m$ modes. 

\subsection{Geometric Programming-based Optimization}
\label{sec:gp}

Our multi-mode optimization problem can be efficiently solved by GP, which is a mathematical optimization method for solving specially formed optimization problems through logarithmic transformations into convex ones. As a result, GP always finds the (true, globally) optimal solution when the problem is feasible~\cite{boyd2007tutorial}. To use GP, the objective function and inequality constraints should be constructed by the special form {\em posynomial}, as in $f(x) = \sum_{k=1}^{K} {c_k x_1^{a_{1k}} x_2^{a_{2k}} \cdots x_n^{a_{nk}} }$, with decision variables $x_i$s, non-negative coefficients $c_k$s, and real-valued exponents $\{a_{11}, \cdots, a_{nK}\}$. Our objective functions and constraints are in posynomial forms except for the idle-time CPU power terms in the rightmost part of $P(\pi)$ in Eq.~\eqref{eq:single_mode} and $P(\hat{\pi})$ in Eq.~\eqref{eq:across_mode}. 

However, it can be removed without affecting optimality as long as $\exists s_i \neq s_{min}$. The optimal solutions in such cases always have 100\% system utilization without any idle time. To prove it intuitively, assume the system utilization $U<100\%$, if we pick a certain task $\tau_i$ and decrease $s_i$, thus increasing $e_i$, until $U$ reaches 100\%, the average power of $\tau_i$ will decrease and the idle power term will disappear, eventually saving more energy than the original configuration. However, when we cannot reduce processor speeds ($\forall i: s_i=s_{min}$), the optimal case would have $U<100\%$ as $p_i$ increases. Such cases only occur when the deadline is extremely long after all $s_i$s are bounded by $s_{min}$. We are not considering such extreme cases in this work.

\section{Safe Mode Change for System Reconfiguration}
\label{sec:mode}

For safe system reconfigurations with dynamic deadlines, the followings should be respected even during mode changes:
\begin{itemize}
    \item {\bf Periodicities.} Every task period before and after the mode change instant should be guaranteed, which can be satisfied by the per-task utilization invariability constraint introduced in Section~\ref{sec:across-modes}.
    \item {\bf Deadlines.} Unfortunately, however, the above periodicities do not guarantee dynamic deadlines, which span across multiple tasks possibly with different modes during a mode change.
\end{itemize}
Thus, this section focuses on developing a safe mode change protocol in terms of end-to-end dynamic deadlines based on the already guaranteed per-task periodicity.

\begin{figure}
\centering
     \begin{subfigure}{0.45\textwidth}
         \centering
         \includegraphics[width=1\textwidth]{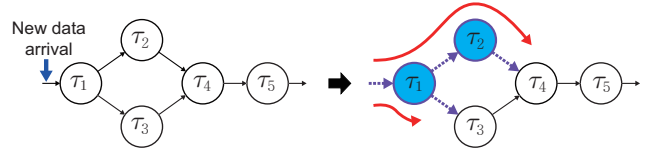}
         \caption{ALAP: tasks gradually change modes (white $\rightarrow$ blue) as the new data progress along paths (curved red arrows) possibly at different speeds.}
         \vspace{0.3\baselineskip}
         \label{subfig-1:alap)}
     \end{subfigure}
     \begin{subfigure}{0.45\textwidth}
         \centering
         \includegraphics[width=1\textwidth]{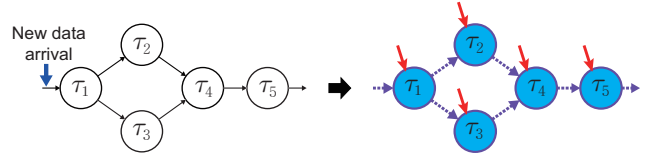}
         \caption{AEAP: the new sensor data arrival immediately triggers (straight red arrows) every task's mode change.}
         \label{subfig-2:aeap}
     \end{subfigure}
        \caption{Baseline mode change protocols}
    \label{fig:mcp}
\vspace{-0.3cm}
\end{figure}

Assume the system is switching from an {\em old} mode to a {\em new} mode, represented by each mode's shortest deadlines, respectively, by
\begin{equation}
\begin{aligned}
d^{old}  &\rightarrow d^{new}.
\end{aligned}
\end{equation}
When $d^{old} < d^{new}$, it is termed as {\em relaxing deadlines} and in the opposite as {\em shrinking deadlines}. System configurations for each mode are denoted by $\pi^{old}\rightarrow\pi^{new}$, that is,
\begin{equation}
\begin{aligned}
((p_1^{old}, s_1^{old}), \cdots, (p_n^{old}, s_n^{old}))&\rightarrow\\
&((p_1^{new}, s_1^{new}), \cdots, (p_n^{new}, s_n^{new}))
\end{aligned}
\end{equation}
in its expansion form. A mode change is triggered by a sensor data arrival at time $t_0$ with its dynamic deadline falling above or below the old range. Then, we consider the following mode change methods, as depicted in Fig.~\ref{fig:mcp}:
\begin{itemize}
    \item {\bf ALAP (As Late As Possible)} individually triggers per-task mode changes after the new data make progress to every incoming edge of it. The actual mode changes will happen at the nearest period boundary after the trigger.
    \item {\bf AEAP (As Early As Possible)} immediately triggers every task regardless of the new data's progress. The mode change completes by $t_0+\text{max}_i(p_i^{old})$ in the worst case when the longest period of $\pi^{old}$ began right before $t_0$.
\end{itemize}
Then we deal with the two cases: (i) relaxing and (ii) shrinking deadlines with the above methods, respectively.

{\bf (i) Relaxing deadline.} In this case, we use {\bf ALAP} such that the mode changes do not adversely affect the already ongoing progress of old sensor data. 
Besides, we need to ensure that the new sensor data do not violate $d^{new}$. 
Note that the new sensor data may progress through tasks possibly with different modes, which happens due to different speeds of different paths. For example, in Fig.~\ref{fig:mcp}(a), $\tau_4$ has two incoming edges, where the upper path requests the mode change while the slower lower path still retains the old mode. 
Then the upper path $\tau_1\rightarrow\tau_2\rightarrow\tau_4\rightarrow\tau_5$ can have a mixture of both modes while handling the new sensor data. 
Thus, the worst-case delay for the new sensor data during {\bf ALAP} mode changes can be calculated as in the following:
\begin{equation}
\begin{aligned}
    D^{new}(\pi^{old} \rightarrow \pi^{new})=&\max_{\forall \delta \in \mathbb{P}}\left(\sum_{i \in \delta}2\max(p_i^{old}, p_i^{new})\right)\\ &\leq\max_{\forall \delta \in \mathbb{P}}\left(\sum_{i \in \delta}2p_i^{new}\right) = d^{new},
\end{aligned}
\label{eq:alap}
\end{equation}
which is less than $d^{new}$ since $\forall i: p_i^{old} \leq p_i^{new}$ that is true when relaxing deadlines. 
This is because, when the task utilization is fixed, the value of $p_i s_i$ should remain the same across modes. Therefore, $p_i$ should increase monotonically to decrease $s_i^{\gamma-1}$ in Eq.~\eqref{eq:power}, reducing the average power in the next longer deadline mode.

{\bf (ii) Shrinking deadlines.} In this case, which basically makes the situation more challenging, we use {\bf AEAP} to quickly finish mode changes, minimizing possible extra delays. Delays for the old sensor data are naturally kept less than $d^{old}$ by the same rationale in Eq.~\eqref{eq:alap} since $\forall i: p_i^{old} \geq p_i^{new}$ when shrinking deadlines. However, regarding the new sensor data, it can suffer extra delays if any task's old period instance that began before the new sensor data arrival persists long enough such that the new data's progress is unexpectedly delayed by that persisting old task instance. Algorithm~\ref{alg:alg} calculates the worst-case delay considering such negative effects for each path $\delta$, which is an ordered set of task indices in each path. Among the calculated delays, we can find the longest. The algorithm gradually accumulates delays by tasks in $\delta$. {\bf Line 1} indicates that it is unavoidable for the new sensor data to be {\em waited} by the old period at the first task. Then we have two cases for the remaining tasks: (i) its mode is already changed before the data progress arrives ({\bf Line 4}) and (ii) an old instance persists ({\bf Line 6}). In the former, we simply accumulate the new delay component $2p_i^{new}$. In the latter, the persisting old (long) period $p_i^{old}$ hides away the accumulated delay up to then, resetting it to $p_i^{old} + p_i^{new}$. 
\begin{algorithm}
\begin{algorithmic}[1]
\REQUIRE $\{(p_1^{old},\cdots, p_n^{old}), (p_1^{new}, \cdots, p_n^{new}), \delta\}$
\ENSURE  The worst-case delay of new sensor data for $\delta$
\STATE $D \leftarrow p_{\delta[1]}^{old} + p_{\delta[1]}^{new}$
\ \ \ \ \ \ \ $\triangleright$ $\delta[1]$ denotes its first element\\
\FOR {$i \in \delta \setminus \{\delta[1]\}$} 
    \IF {$D > p_i^{old} - p_i^{new}$} 
        \STATE $D \leftarrow D + 2p_i^{new}$
    \ELSE 
        \STATE $D \leftarrow p_i^{old} + p_i^{new}$
    \ENDIF
\ENDFOR
\RETURN$D$
\end{algorithmic}
\caption{Finding the worst-case delay for {\bf AEAP}}
\label{alg:alg}
\end{algorithm}

By the above analyses, we claim that when relaxing deadlines with {\bf ALAP}, there is no deadline miss for both the already ongoing progress and new ones. When shrinking deadlines with {\bf AEAP}, the already ongoing progress rather benefits from it, whereas new sensor data can suffer extra delays. However, we can analyze the worst-case extra delays, which can be used when planning appropriate safety margins in design time.

\begin{figure*}
  \includegraphics[width=\textwidth,height=4cm]{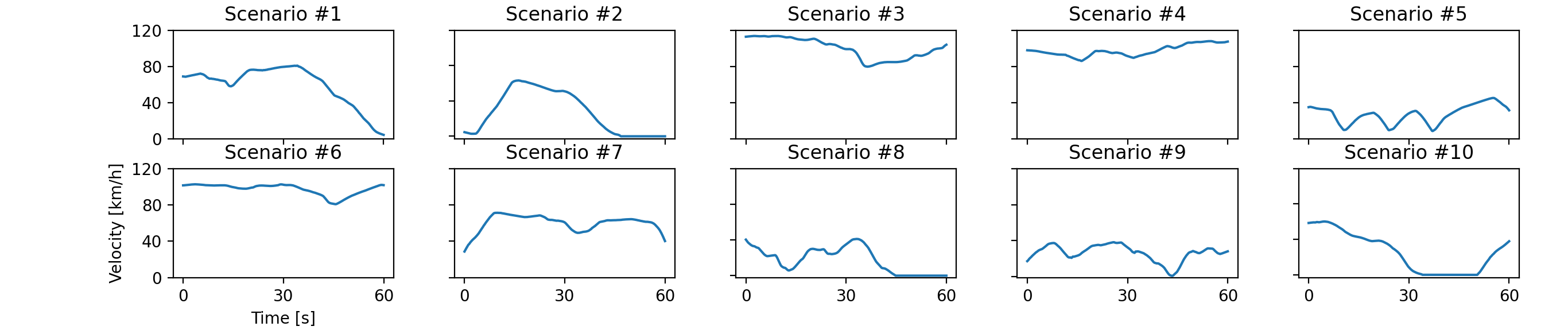}
  \caption{Real-world driving scenarios from comma.ai dataset}
\label{fig:scenarios}
\end{figure*}

\section{Practical Considerations}
\label{sec:practical}

\subsection{Extension to Discrete CPU Frequency Levels}
\label{sec:practical-discrete}

Our optimization method yields speed factors in the continuous range between $s_{min}$ and 1. However, because most CPUs, in practice, support only a predefined set of discrete frequency levels, we need to adapt the resulting speed factors to the discrete domain. One possible approach is to emulate the exact speed factors by modulating between two neighboring discrete frequency levels~\cite{yasuura_intradvfs_1998, bini_intradvfs_2005} to obtain the near-optimal energy reduction similar to the continuous frequency solution. However, using such intra-task DVFS entangles other practical considerations such as extra time and energy overhead associated with excessive frequency transitions~\cite{kim_dvfs_2004} and possible transient faults~\cite{zhang_intradvfs_2003, zhu_intradvfs_2004}.

In light of this, we propose a more conservative but safer method that uses the closest frequency level that is higher than the corresponding optimal speed factor. Then every task's actual utilization is less than or equal to the ideal utilization in accordance with the per-task utilization invariability constraint. Thus, even though the approximated system will consume more energy than the continuous one, it safely ensures schedulability and end-to-end delay constraint during mode changes.



\subsection{Extension to Multicore CPUs}

Our system model assumes a singlecore CPU, which is still dominant in the automotive industry due to the safety concerns and software complexity caused by multicore CPUs. Nonetheless, high performance multicore CPUs have recently gained acceptance in the industry, necessitating the extension of our method to multicore CPUs. Although this issue is not within this paper's scope, this section briefly discusses some ideas for the extension.

One possible approach is to incorporate one of the global multicore scheduling algorithms that still provide utilization-based schedulability analysis. Then, only minor modifications are necessary in the optimization formulation, and also the mode change protocols can be easily reused. However, this extension approach requires multiple tasks executing in parallel to have different CPU speed levels, which is difficult to implement because most CPUs do not support per-core DVFS in practice.

Another approach, which seems more practical, is to use the RT-Gang scheduling framework~\cite{ali_rt-gang_2019} into our method. RT-Gang may eliminate many multicore issues including per-core DVFS and unpredictable inter-core memory bandwidth contention by its one-gang-at-a-time policy. Here, a gang is a predefined group of tasks that is allowed to run in parallel. After the gang grouping, the system model is transformed as in our singlecore model, where despite the gang grouping problem, our optimization method and mode change protocol can be used without modification. For the remaining gang grouping, we leave it for future work.



\section{Experiments}
\label{sec:experiments}

\begin{table*}[ht]
\caption{Workload information at maximum speed ($s_i=1$)}
\centering
\begin{tabular}{c||c|c|c|c|c|c|c|c|c|c}
\hline
{\bf Task} &  Camera Grabber & Lidar Grabber & CAN & SFM & Lane detection & Detection & Loc & EKF & Planner & DASM \\
\hline
{\bf WCET (ms)} & 0.25 & 2.75 & 0.15 & 6.95 & 10.55 & 29.00 & 73.70 & 1.10 & 3.10 & 0.325 \\
\hline
\end{tabular}
\label{tab:task}
\vspace{-0.5cm}
\end{table*}

\subsection{Experimental Setup}
\label{sec:expr_setup}
{\bf Workload.} We use the DAG in Fig.~\ref{fig:Dag} with ten tasks, where their WCETs are listed in Table.~\ref{tab:task}. They are derived from the actual measurements~\cite{krawczyk_analytical_2019} on Nvidia Jetson TX2 platform (Denver cores) and scaled considering typical high-performance computing systems for autonomous driving. 
Unfortunately, these industry-level applications are IP protected, so we were not able to run them on a real platform.
\\
{\bf Energy Model.} We empirically found the energy parameters for Eq.~\eqref{eq:power}, as $\alpha$=842.04, $\beta$=232.81, and $\gamma$=2.64, on the same hardware platform in~\cite{krawczyk_analytical_2019}.\\
{\bf Discrete frequency levels.} We use 12 evenly spaced frequencies between 345~MHz and 2~GHz from the same hardware platform in~\cite{krawczyk_analytical_2019}.\\
{\bf Scenarios.} We use real-world driving scenarios from the comma.ai driving dataset~\cite{velocity_dataset}, where we picked ten 60-second driving logs with their velocity from 0~km/h to 114~km/h as depicted in Fig.~\ref{fig:scenarios}.  \\
{\bf Dynamic Deadlines.} We converted velocities into deadlines using Eq.~\eqref{eq:dv} with $\lambda$=20~m and $a^{max}$=2.5~ms\textsuperscript{-2}. The shortest deadline (617~ms) is obtained using the highest velocity (114~km/h), while the longest deadline (2945~ms) is bounded by $s_{min}$=0.17. Then, the number of modes is chosen arbitrarily, $m$=24, partitioning the range with equal length. \\
{\bf Optimization.} For the GP solver, we use CVX~\cite{cvx} on a laptop with a quad-core Intel Core i7@2.6 GHz CPU, where our optimization takes about 14.6 seconds.\\
{\bf Simulation.} We implemented a simulator supporting the EDF scheduling and our mode change protocol, by which we can precisely estimate the exact task schedules, end-to-end delays, and energy consumption.


\subsection{Evaluation}
\label{sec:eval}

\begin{figure}
\centering
\includegraphics[width=9.0cm]{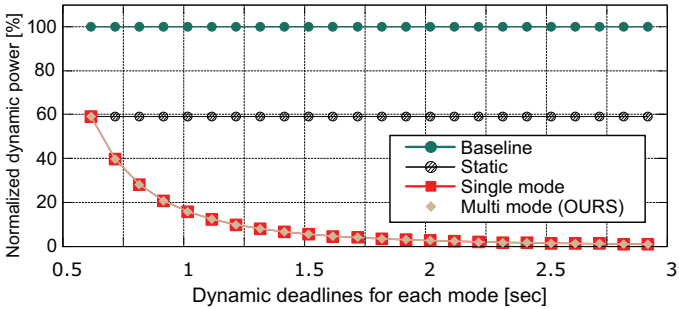}
\caption{Dynamic power optimization results}
\label{fig:single_vs_across}
\vspace{-0.4cm}
\end{figure}

The goal of our evaluation is twofold: 
(i) illustrate the benefit of using multiple frequency levels in dynamic deadline situation over conventional methods that utilize only one or two frequencies 
(ii) evaluate the effectiveness of our mode change protocol. It is of no use if it compromises safety. The following three methods are compared in the evaluation:
\begin{itemize}
\item {\bf Baseline}: $\forall i: s_i=1$; 
\item {\bf Static}: optimized for the minimum deadline (617~ms);
\item {\bf Multi mode (OURS):} method in Subsection~\ref{sec:across-modes}.
\end{itemize}

Fig.~\ref{fig:single_vs_across} shows the dynamic power optimization results of each method. 
For \textit{Single mode}, we solved the problem in Subsection~\ref{sec:per-mode} repeatedly for each deadline.
In the leftmost mode with the shortest deadline, the three methods except the baseline method show an equal result (58.9\%). 
However, as the deadline increases, the dynamic power plunges to 1.0\% in the rightmost mode with the longest deadline. 
Another interesting observation is that our multi-mode results do not reveal visible performance degradation compared with the single-mode results, meaning that the per-task utilization invariability constraint does not critically affect the optimization results. 
On top of that, there is a 75\% reduction in overall optimization time when using multi-mode formulation thanks to fewer decision variables.

\begin{figure}
\centering
\includegraphics[width=9cm]{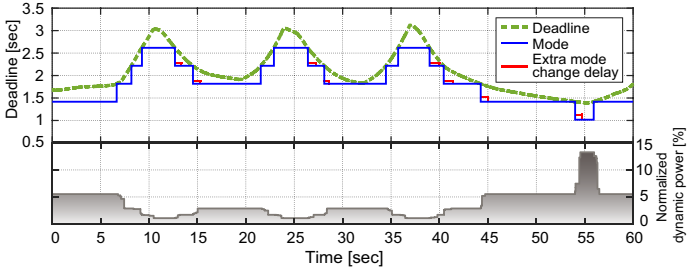}
\caption{Simulation results with a real-world driving scenario}
\label{fig:driving_scenario_example}
\vspace{-0.4cm}
\end{figure}

Fig.~\ref{fig:driving_scenario_example} shows the simulation results with one of our driving scenarios. Only in this experiment, we partition the deadline range into six equal length modes for better visualization. As the deadline changes, as depicted by the green dashed line, mode changes are sporadically triggered, depicted by the rising and falling edges of the blue line that roughly follows the dynamic deadline. Note the small red staircase shapes at each falling edge, which depicts the extra delays when shrinking deadlines analyzed by Algorithm~\ref{alg:alg}. More specifically, their height represents the extra delay, while their width represents the transient interval for each mode change. As shown in the figure, we do not have any deadline miss in this scenario even after considering the extra mode change delays. Additionally, the bottom plot shows how the normalized dynamic power varies (between 1.0\% and 13.3\%) depending on the mode changes.


\begin{figure}
\centering
\includegraphics[width=9cm]{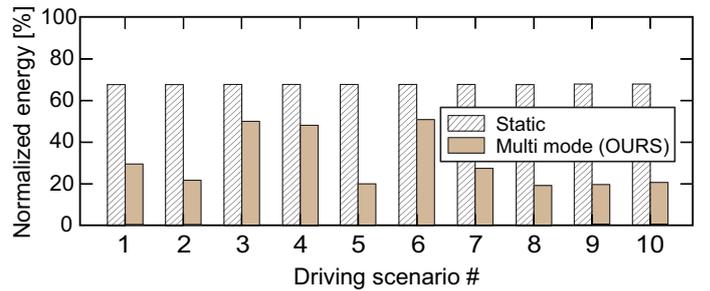}
\caption{Energy consumption with varying driving scenarios (continuous frequency)}
\label{fig:driving_scenario_compare_cont}
\vspace{-0.4cm}
\end{figure}

\begin{figure}
\centering
\includegraphics[width=9cm]{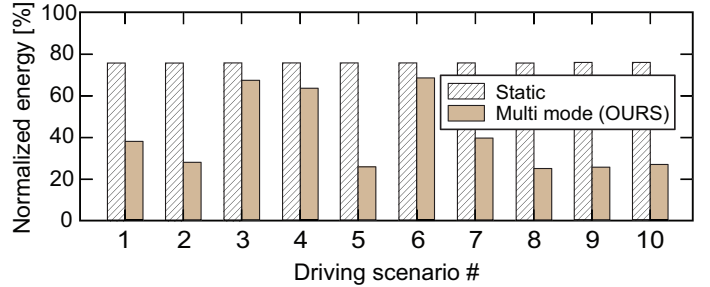}
\caption{Energy consumption with varying driving scenarios (discrete frequencies)}
\label{fig:driving_scenario_compare_disc}
\vspace{-0.4cm}
\end{figure}

Fig.~\ref{fig:driving_scenario_compare_cont} shows the energy consumption, including static and dynamic power, in driving scenarios using continuous frequency.
Though \textit{Static} achieves a significant energy reduction of 32\% from baseline, \textit{Multi mode} reduces further by timely utilizing lower frequencies. 
Scenarios \#3, \#4, and \#6 show relatively low energy reductions since they maintain the vehicle within a high-velocity range most of the time, making the system remain in the short deadline modes. In the others, more than half of the total energy was saved. On average, our method achieved 69.3\% and 54.8\% energy reductions compared with the baseline and the static methods, respectively. 

Fig.~\ref{fig:driving_scenario_compare_disc} shows the energy consumption using discrete frequencies as discussed in Section~\ref{sec:practical-discrete}.
A discrete frequency that is greater than or equal to the optimal one is assigned to each task. 
Therefore, we expect the energy consumption will increase in both methods.
\textit{Static} achieves 24\% of energy reduction, which is smaller than the previous case. \textit{Multi mode} still shows a significant reduction from \textit{Static}, however, the energy gap between them becomes smaller. In scenarios \#3, \#4, and \#6, the reduced gap is very notable. 
This is caused by the nonlinearity in the power model in terms of speed factor.
The difference between the optimal and mapped discrete frequency yields bigger difference in power consumption in the higher frequency range. 
On average, there was a 31.8\% increase in energy consumption when mapped to discrete frequencies.
Our method achieved 59.3\% and 46.4\% energy reductions on average compared with the baseline and the static methods, respectively. 




Throughout the experiments with the industry-level applications and real-life driving scenarios, we could not find any deadline violation with our mode change protocols. 
However, as we discussed in Section ~\ref{sec:mode}, the extra delays when shrinking deadlines are unavoidable and they could be an issue depending on vehicle's maximum acceleration and mode lengths. In that case, the violations can be avoided by selecting a shorter deadline mode considering possible mode change delays analyzed in Section \ref{sec:mode} (i.e. a safety margin) in design time analysis. 
We leave the choice of optimal margin for our future work since it requires a thorough analysis.

\section{Related Work}
\label{sec:related}

{\bf Dynamic deadlines.}
Recent studies~\cite{lee2020subflow, heo2020real} presented motivating examples of dynamic deadlines in autonomous driving, where object detection systems are commonly proposed that adapt themselves to varying deadlines, demonstrating the unique potential of autonomous driving systems. More specifically, Lee and Nirjon~\cite{lee2020subflow} support dynamic deadlines with selective subgraph executions by considering varying time budgets. Heo et al.~\cite{heo2020real} support dynamic deadlines by selectively executing multiple forward propagation paths of a neural network with different execution times. Both studies trade dynamic deadlines (or slacks) for improving the object detection accuracy. However, little work has been done with dynamic deadlines for the energy optimization of autonomous driving.

{\bf DVFS-based energy optimization.}
There have been many efforts to develop energy-efficient real-time systems, most of which, however, assume only static deadlines. Broadly, there are two frequently used energy saving approaches in hard real-time systems: DVFS and {\em dynamic power management} (DPM). In DVFS approaches, there is a body of literature to find speed levels through offline optimization \cite{daren_time_2014, bambagini2016energy, guo_energy-efficient_2019, saifullah2020cpu}.
These approaches try to find critical speed factors under a static deadline constraint, which is the lowest frequency satisfying the given static timing constraint. In contrast, our method finds the critical speed for each deadline through multiple modes with different timing constraints and employs a safe mode change protocol to freely go back and forth between them. Another body of real-time DVFS approaches tries to reclaim \textit{dynamic slacks}\cite{jejurikar2005dynamic, bambagini_dvfs_2016}, which is not to be confused with our dynamic deadlines; instead, they define dynamic slacks as the difference between the worst-case and the actual execution times. Note that dynamic slack reclaiming does not conflict with our approach and could be used together for further energy reductions.

{\bf DPM-based energy optimization.} In DPM approaches, cores are switched off during idle periods to reduce energy consumption.
It is useful when the system is underutilized. 
However, each idle period should be large enough to offset the energy overhead of frequent switching on and offs.  
Many scheduling methods have been proposed to create large idle periods \cite{yann-hang_lee_scheduling_2003, huang_periodic_2009, awan_enhanced_2011, sun_flow_2019}. 
Lee et al.\cite{yann-hang_lee_scheduling_2003} proposed leakage control EDF scheduling, but it requires additional hardware to calculate the sleep period. 
This impractical assumption was avoided in \cite{awan_enhanced_2011} by using a simpler sleep period calculation method and the enhanced race-to-halt (ERTH) algorithm. 
With ERTH, tasks run at full speed to secure longer sleep time.
DPM and DVFS approaches are not orthogonal and can be used in conjunction.
However, satisfying the timing constraints during mode changes with DPM approaches requires non-trivial considerations.
Thus, we leave integrating DPM for our future work.




\section{Conclusion}
\label{sec:conclusion}

This study is motivated by emerging autonomous driving applications with time-varying dynamic deadlines, where the computing system's excessive energy consumption is a major concern. Unlike traditional energy optimization methods assuming rigid static deadlines, our solution approach tries to utilize the dynamic slack obtained by adaptively relaxing deadlines considering the vehicle's physical state. For that, our GP-based optimization method proactively exploits the dynamic deadlines to find energy-efficient multi-mode system configurations. Moreover, our safe mode change protocol enables adaptive system reconfiguration between the predefined modes. Our experimental results show an average of 46.4\% energy reduction from the previous method with static deadlines, demonstrating the great potential toward energy-efficient autonomous driving systems.

Based on the theoretical foundation presented in this paper, in the future, we plan to extend our work to incorporate other energy saving techniques into our adaptive reconfiguration and to consider practical issues in COTS CPUs including multicore processors. Also, we plan to extend our system model considering other accelerators such as GPUs.




\section*{Acknowledgment}
This work was supported partially by NSF grant CCF-1704859 and partially by the Ministry of Land, Infrastructure, and Transport (MOLIT), Korea, through the Transportation Logistics Development Program (20TLRP-B147674-03, Development of Operation Technology for V2X Truck Platooning). J.-C. Kim is the corresponding author.
\balance
\bibliographystyle{IEEEtran}
\bibliography{energy}

\end{document}